\documentclass[11pt]{amsart}
\usepackage{amsfonts}
\usepackage{amsmath}
\usepackage{amssymb}
\usepackage{amsthm,amsfonts,latexsym,epsfig,geometry}
\usepackage{hyperref,times}
\usepackage{amscd}
\usepackage{mathtools}
\usepackage{color}

\allowdisplaybreaks[1]
\newtheorem{theorem}{Theorem}[section]
\newtheorem{lemma}[theorem]{Lemma}

\theoremstyle{definition}

\theoremstyle{remark}

\numberwithin{equation}{section}
\usepackage{hyperref}
\hypersetup{ 
	colorlinks=true,
	linkcolor={blue}, urlcolor={blue},
	citecolor={red}, anchorcolor = {blue}}
\begin{document}
\title[ ]{On Some New Congruences For Biregular Overpartitions}
\author{ANAKHA V }
\address{ANAKHA V, NATIONAL INSTITUTE OF TECHNOLOGY CALICUT, KOZHIKODE, KERALA, INDIA 673601 }
\email{anakha$\textunderscore$p240019ma@nitc.ac.in }
 
\begin{abstract}

Inspired by the recent work by Nadji, Ahmia and Ramírez \cite{biregularNadji}, we examined the arithmetic properties of $\overline{B}_{\ell_1,\ell_2}(n)$, the number of overpartitions of $n$ whose parts are neither divisible by $\ell_1$ nor divisible by $\ell_2$. In particular, we establish some congruences modulo $k\in\{4,8,6,12\}$ satisfied by $\overline{B}_{\ell_1,\ell_2}(n)$ where $\ell_1$ and $\ell_2$ take values as arbitrary powers of $2$ and $3$. Moreover, we extend certain results proved in \cite{biregularNadji} and \cite{ghoshal2025arxiv} for $\ell_1$ and $\ell_2$ with random powers of $2$ and $3$. Generating functions, dissection formulas, and theta functions are used to prove our main findings.

\end{abstract}
\maketitle
\section{Introduction and statement of results}
\label{intro}
 Suppose that $n$ is a positive integer. A \emph{partition} of $n$ is a non-increasing sequence of positive integers $(\lambda_1,\lambda_2,\hdots,\lambda_k)$ such that $n=\lambda_1+\lambda_2+\cdots+\lambda_k$. Each of the $\lambda_i$ are called a \emph{part}. Let the number of partitions of $n$ is denoted by $p(n)$. By convention, $p(0)=1$. The generating function for $p(n)$ is given by
 \begin{align*}
    \sum_{n=0}^{\infty}{p}(n)q^n=\prod_{n=1}^\infty\frac{1}{1-q^n}=\frac{1}{(q;q)_\infty}.
\end{align*}

Here and throughout this paper we will use the standard $q-$series notation for $q\in\mathbb{C}$, $|q|<1$:

\begin{equation*}
\begin{split}
    (a;q)_0&:= 1 \\
     (a;q)_n&:=   \prod_{k=0}^{n-1}(1-aq^k) \text{ for }n\geq 1 \text{ and }\\
    (a;q)_\infty&:=\prod_{k=0}^{\infty}(1-aq^k).
\end{split} 
\end{equation*}

 An \emph{overpartition} of $n$ is a partition of $n$ in which the first occurrence of each part may be overlined. The number of overpartitions of $n$ is denoted by $\overline{p}(n)$, and by convention $\overline{p}(0)=1$. As an example, $\overline{p}(4)=14$ and the overpartitions are 
\begin{align*}
    (4),(\overline{4}), (3,1),(\overline{3},1),(3,\overline{1}),(\overline{3},\overline{1}),(2,2), (\overline{2},2), (2,1,1),(\overline{2},1,1), (2,\overline{1},1), (\overline{2},\overline{1},1), (1,1,1,1), (\overline{1},1,1,1).
\end{align*}

The generating function for the number of overpartitions is given by Corteel and Lovejoy \cite{corteel2004overpartitions},
\begin{align*}
    \sum_{n=0}^{\infty}\overline{p}(n)q^n=\prod_{n=1}^\infty\left(\frac{1+q^n}{1-q^n}\right)=\frac{(q^2;q^2)_\infty}{(q;q)_\infty^2}.
\end{align*}

Extensive studies have been conducted on arithmetic properties of overpartitions. Hirschhorn and Sellers \cite{hirschhorn2005overpartition} proved the following

\begin{align*}
    \overline{p}(9^\alpha(27n+18))\equiv 0\pmod 3
\end{align*}
and

\begin{align*}
    \overline{p}(27n)\equiv \overline{p}(3n) \pmod 3.
\end{align*}

For further reading on overpartitions we cite \cite{kim2008overpartition,kim2009short,xia2017overpartitions,yao2013overpartitions}.\\

For an integer $\ell>1$, an $\ell-$regular partition of $n$ is a partition that has no parts that are divisible by $\ell$. Let $b_\ell(n)$ denote the number of $\ell-$regular partitions of $n$. Then the generating function is given by
\begin{align*}
    \sum_{n=0}^\infty b_\ell(n)q^n=\frac{(q^\ell;q^\ell)_\infty}{(q;q)_\infty}.
\end{align*}
Numerous studies have explored the arithmetic properties of regular partitions. Cui and Gu \cite{cui2013regular} discovered several infinite families of congruences modulo $2$ for $\ell-$regular partition functions for $\ell=2,4,5,8,13,16$. Hirchhorn and Sellers yielded many Ramanujan type congruences for $b_5(n)$ in their work \cite{hirschhorn5regular}.\\

A bipartition of $n$ is a pair $(\lambda,\mu)$ of partitions such that the sum of all its parts is equal to $n$. In \cite{linl-regularbipartition}, Lin studied the properties of $B_{\ell}(n)$ which is the number of $\ell-$regular bipartitions of $n$ and having the generating function
\begin{align*}
    \sum_{n=0}^\infty B_\ell(n)q^n=\frac{(q^\ell;q^\ell)_\infty^2}{(q;q)_\infty^2}.
\end{align*}
In particular, he proved
\begin{align*}
    B_7\left(3^\alpha n+\frac{5\cdot3^{\alpha-1}-1}{2}\right)\equiv 0\pmod 3.
\end{align*}
Let $B_{\ell,m}(n)$ denotes the count for $(\ell,m)-$regular bipartitions of $n$ which is defined as a bipartition $(\lambda,\mu)$ of $n$ such that $\lambda$ is $\ell-$regular and $\mu$ is $m-$regular with a generating function given by

\begin{align*}
     \sum_{n=0}^\infty B_{\ell,m}(n)q^n=\frac{(q^\ell;q^\ell)_\infty (q^m;q^m)_\infty}{(q;q)_\infty^2}.
\end{align*}
In \cite{dou2016lm-regularbipartition} Dou proved 
\begin{align*}
    B_{3,11}\left(3^\alpha n+\frac{5\cdot3^{\alpha-1}-1}{2}\right)\equiv 0\pmod {11}.
\end{align*}
Some more arithmetic propertities of $B_{\ell,m}(n)$ can be found in \cite{wang2017arithmetic,adiga2019st-regularbipartition}.

Lovejoy \cite{lovejoy2003gordon} investigated the $\ell-$regular overpartition function $\overline{A}_\ell(n).$ Shen \cite{shen2016arithmetic}
 provided the generating function for $\overline{A}_\ell(n)$ as 
\begin{align*}
    \sum_{n=0}^\infty \overline{A}_\ell(n)q^n=\frac{(q^\ell;q^\ell)_\infty^2(q^2;q^2)_\infty}{(q;q)_\infty^2(q^{2\ell};q^{2\ell})_\infty}.
\end{align*}
and derived some congruences for $\overline{A}_3(n)$ and $\overline{A}_4(n)$. Saikia and Barua \cite{saikia2017congruences} proved congruences for $\overline{A}_5(n)$ modulo $4,$ $\overline{A}_6(n)$ modulo $3$, and $\overline{A}_8(n)$ modulo $4.$ One can refer \cite{shenlregular,alanazi2016infinite}, for more information. \\
 Andrews \cite{andrews2015singular} introduced singular overpartition function $\overline{C}_{k,i}(n)$ which enumerates the number of overpartitions of $n$ in which the parts are not divisible by $k$ and the only overlined parts are those which are congruent to $\pm i$ modulo $k$. The generating function for $\overline{C}_{k,i}(n)$ is given as 
 \begin{align*}
     \sum_{n=0}^{\infty}\overline{C}_{k,i}(n)q^n=\frac{(q^k;q^k)_{\infty}(-q^i;q^k)_{\infty}(-q^{k-i};q^k)_{\infty}}{(q;q)_{\infty}}
 \end{align*}
 for $k\geq 3$ and $1\leq i\leq \lfloor \frac{k}{2}\rfloor$.
 In the same article Andrews proved the following congruences:
 \begin{align*}
   \overline{C}_{3,1}(9n+3)\equiv \overline{C}_{3,1}(9n+6)\equiv 0 \pmod 3 \text{ \hspace{1cm} for all $n\geq 0$ }.  
 \end{align*}
 Singular overpartitions is a well established area of research. Chen, Hirschhorn and Sellers \cite{chen2015singular} discovered congruences for $\overline{C}_{3,1}(n)$, $\overline{C}_{4,1}(n),$ $\overline{C}_{6,1}(n)$ and $\overline{C}_{6,2}(n)$ modulo $3$ and powers of $2$. Later, Ahmed and Baruah \cite{ahmed2015singular} presented certain congruences for $\overline{C}_{3,1}(n),\overline{C}_{8,2}(n),\overline{C}_{12,2}(n),\overline{C}_{12,4}(n),\overline{C}_{24,8}(n) $ and $\overline{C}_{48,16}(n)$. In \cite{ray2024distribution}, Ray studied the distribution for $\overline{C}_{p,1}(n),$ where $p\geq 5$ is a prime number. The reader can look into \cite{barmanRay2019divisibility,singhBarman2021new} for more insights.
 
 Bringmann \emph{et al.} \cite{bringmann2015overpartitions} proposed and studied a new type of overpartitions characterized by restricted odd difference. It further researched by Hirschhorn and Sellers \cite{hirschhorn2020restricted} and proved some parity results. Many other results on restricted overpartitions can be found in the literature. Referrence \cite{baruah2025restricetd} provide a relevant  example examining the arithmetic properties of restricted overpartitions, specially focussing on parts belonging to certain residue classes modulo $8$. 

Recently, Nadji, Ahmia, and Ramrez \cite{biregularNadji} have delved into the arithmetic properties of biregular overpartitions, which is defined as follows:\\
For a pair of two relatively prime integers $(\ell_1,\ell_2)>1$, an $(\ell_1,\ell_2)-$\emph{biregular overpartitions} of $n$ is a partition in which none of the parts are divisible by $\ell_1$ or $\ell_2$.\\
They provided some congruences modulo $3$ and powers of $2$ for the values of pairs \hspace{0.06cm} $(\ell_1,\ell_2)\in \{(4,3),(4,9),\\
(8,3),(8,9)\}$. For proving these congruences they used generating functions, dissection formulas and Smoot's implementation of Radus's Ramanujan-Kolberg algorithm \cite{smoot2021computation}. To know more about Radu's algorithm, the reader can take a look into \cite{radu2015algorithmic1}. In \cite{manjil2024some}, Alanazi, Munagi and Saikia provided numerous congruences for the biregular overpartitions, which they showed through modular forms and generating functions. They also presented combinatorial proof for the congruences modulo $4$. They proved both analytically and combinatorially the following:
\begin{align*}
    \overline{B}_{\ell_1,\ell_2}(n)\equiv 0 \pmod{2}\text{ \hspace{1cm} for $n\geq 1$}.
\end{align*}
Subsequently, Ghoshal and Jana \cite{ghoshal2025arxiv} investigated additional properties of biregular overpartitions. To elaborate, through elementary methods they derived some of the congruences given in \cite{biregularNadji} in which the proof is done by making use of Radu's algorithm. Moreover, they developed a general method for proving congruences modulo $8$ which is applicable without any specific $2-$dissection requirements. Paudel, Sellers and Wang \cite{paudel2025extending} extended several results of Alanazi \emph{et al.} and established a large number of new congruences. They proved their findings  merely by making use of some $q-$series manipulations and dissection formulas.In \cite{newbiregular}, the authors have derived congruences for $(\ell_1,\ell_2)\in \{(2,9),(5,2),(8,3)\}$ and generally for $(\ell_1,\ell_2)\in \{(5,2^t),(3,2^t),(4,3^t)\}$.

In our paper, we follow the same notation $\overline{B}_{\ell_1,\ell_2}(n)$ as \cite{biregularNadji} for the number of $(\ell_1,\ell_2)-$biregular overpartitions of $n$.

The generating function for $\overline{B}_{\ell_1,\ell_2}(n)$ is 
\begin{equation}
\label{0}
    \begin{split}
        \sum_{n=0}^{\infty}\overline{B}_{\ell_1,\ell_2}(n)q^n&=\prod_{n=1}^\infty \frac{(1+q^n)(1+q^{\ell_1\ell_2 n})(1-q^{\ell_1 n})(1-q^{\ell_2 n})}{(1+q^{\ell_1 n})(1+q^{\ell_2 n})(1-q^n)(1-q^{\ell_1\ell_2 n})}\\
        &=\frac{(q^2;q^2)_\infty(q^{\ell_1};q^{\ell_1})_\infty^2(q^{\ell_2};q^{\ell_2})_\infty^2(q^{2\ell_1\ell_2};q^{2\ell_1\ell_2})_\infty}{(q;q)_\infty^2(q^{2\ell_1};q^{2\ell_1})_\infty(q^{2\ell_2};q^{2\ell_2})_\infty(q^{\ell_1\ell_2};q^{\ell_1\ell_2})_\infty^2}.
    \end{split}
\end{equation}

We give a certain kind of generalization for some of the results proved in \cite{biregularNadji} and \cite{ghoshal2025arxiv}. Besides that our new findings are also included. To be more exact, we present the proof of the following results in which we make use of generating functions and dissection formulas.

\begin{theorem}\label{thm1}
    For every integers $n\geq 0$ and for $\alpha\geq 2$, we have
    \begin{align}
    \label{7}
        \overline{B}_{2^\alpha,3}(4n+2)\equiv 0\pmod 4,\\
    \label{8}    
        \overline{B}_{2^\alpha,3}(8n+5)\equiv 0\pmod 4,
        \end{align}
 \hspace{1cm}  and for $\alpha\geq 3$
     \begin{align}  
     \label{9}   
        \overline{B}_{2^\alpha,3}(8n+6)\equiv 0\pmod 8.
    \end{align}
\end{theorem}

  \begin{theorem}\label{thm2}
  
      For every integer $n\geq 0$ and $\alpha>2k +1$, where $k\geq 0$ we have
\begin{align}\label{th2e1}
\overline{B}_{2^\alpha,3}(4^k(4n+3))\equiv 0 \pmod 6
\end{align}
 \hspace{1cm}  For $\alpha>2,$
     \begin{align} 
\label{th2e2}
\overline{B}_{2^\alpha,3}(8n+7)\equiv 0 \pmod {12}.
\end{align}

  \end{theorem}

  \begin{theorem}\label{thm3}
      Let $n\geq 0$  and $\beta\geq 2$ be integers. Then for any integer $\alpha \geq 0$, we have
      \begin{align}\label{th3e1}
       \overline{B}_{2^{2\alpha +1},3^\beta}(9n+3i)\equiv 0 \pmod 4  
      \end{align}
      and for $\alpha \geq 1$, we have
\begin{align}\label{th3e2}
    \overline{B}_{2^{2\alpha},3^\beta}(9n+3i)\equiv 0 \pmod 8
\end{align}
   for $i=1,2.$   
  \end{theorem}

We make use of the proof techniques given in \cite{ghoshal2025arxiv} to prove the following theorem.

\begin{theorem}\label{THM}
    Let $\alpha\geq 2$, $\beta\geq 1$ and $n\geq 0$ be integers. If $\alpha$ is odd and $\beta$ is even or both $\alpha$ and $\beta$ are even, then
     \begin{align}\label{THe4}
        \overline{B}_{2^\alpha,3^\beta}(12n+3)\equiv0 \pmod 8\\
         \label{THe1}  \overline{B}_{2^\alpha,3^\beta}(12n+7)\equiv0 \pmod 8.
    \end{align}
    If $\alpha$ and $\beta$ have different parity or both $\alpha$ and $\beta$ are even, then
     \begin{align}
\label{THe2} \overline{B}_{2^\alpha,3^\beta}(12n+11)\equiv0 \pmod 8.
    \end{align}
    
\end{theorem}

\section{Preliminaries}

In this section, we state some lemmas and formulas which we use to prove our results. Ramanujan's theta function \cite[pp. 34, 18.1]{ramanujan3} $f(a,b)$ is defined by
\begin{align*}
    f(a,b)=\sum_{n=-\infty}^{\infty}a^{n(n+1)/2}b^{n(n-1)/2}=(-a;ab)_\infty(-b;ab)_\infty(ab;ab)_\infty, \text{ for }|ab|<1.
\end{align*}
 One of the important representations of theta function is given as:
 \begin{align*} 
     \varphi(q):=f(q,q)=1+2\sum_{n=1}^{\infty}q^{n^2}=\frac{(-q;q^2)_\infty(q^2;q^2)_\infty}{(q;q^2)_\infty(-q^2;q^2)_\infty}=\frac{(q^2;q^2)_\infty^5}{(q;q)_\infty^2 (q^4;q^4)_\infty^2}.
 \end{align*}
Replacing $q$ by $-q$ in the above equation yields

\begin{align*}
    \varphi(-q)=\frac{(q;q)_\infty^2}{(q^2;q^2)_\infty}.
\end{align*}

\begin{lemma}
    For all primes $p$ and all $k,m\geq 1$, we have
    \begin{align}\label{lp1}
        f_{pm}^{p^{k-1}}\equiv f_m^{p^k}\pmod{p^k}
    \end{align}
\end{lemma}

We list some dissection formulas which we use in our proofs.

\begin{lemma}
    We have
    \begin{align}
\label{1}
        \frac{(q^3;q^3)_\infty^2}{(q;q)_\infty^2}&=\frac{(q^4;q^4)_\infty^4(q^6;q^6)_\infty(q^{12};q^{12})_\infty^2}{(q^2;q^2)_\infty^5(q^8;q^8)_\infty(q^{24};q^{24})_\infty}+2q\frac{(q^4;q^4)_\infty(q^6;q^6)_\infty^2(q^8;q^8)_\infty(q^{24};q^{24})_\infty}{(q^2;q^2)_\infty^4(q^{12};q^{12})_\infty}\\
\label{2}
        \frac{(q^3;q^3)_\infty}{(q;q)_\infty^3}&=\frac{(q^4;q^4)_\infty^6(q^6;q^6)_\infty^3}{(q^2;q^2)_\infty^9(q^{12};q^{12})_\infty^2}+3q\frac{(q^4;q^4)_\infty^2(q^6;q^6)_\infty(q^{12};q^{12})_\infty^2}{(q^2;q^2)_\infty^7}   
    \end{align}
\end{lemma}
Xia and Yao \cite{xia2013analogues} proved equation (\ref{1}). Barua and Ojah \cite{baruah2012analogues} proved (\ref{2}).
\begin{lemma}
    We have
    \begin{align}
\label{3}
       \frac{1}{(q;q)_\infty^4}&=\frac{(q^4;q^4)_\infty^{14}}{(q^2;q^2)_\infty^{14}(q^8;q^8)_\infty^4}+4q\frac{(q^4;q^4)_\infty^2(q^8;q^8)_\infty^4}{(q^2;q^2)_\infty^{10}} \\
\label{4}      
       \frac{1}{(q;q)_\infty^2}&=\frac{(q^8;q^8)_\infty^5 }{(q^2;q^2)_\infty ^5 (q^{16};q^{16})_\infty^2}+2q \frac{(q^4;q^4)_\infty ^2  (q^{16};q^{16})_\infty^2}{(q^2;q^2)_\infty ^5 (q^{8};q^8)}      
    \end{align}
\end{lemma}
Equations (\ref{3}) and (\ref{4}) are consequences of dissection formulas of Ramanujan \cite{ramanujan3}. 

\begin{lemma}
    We have
    \begin{align}
    \label{5}
        \frac{(q;q)_\infty^2}{(q^2;q^2)_\infty}&=\frac{(q^9;q^9)_\infty^2}{(q^{18};q^{18})_\infty}-2q\frac{(q^3;q^3)_\infty(q^{18};q^{18})_\infty^2}{(q^6;q^6)_\infty(q^9;q^9)_\infty}\\
    \label{6}
         \frac{(q^2;q^2)_\infty}{(q;q)_\infty^2}&=\frac{(q^6;q^6)_\infty^4(q^9;q^9)_\infty^6}{(q^3;q^3)_\infty^8(q^{18};q^{18})_\infty^3}+2q\frac{(q^6;q^6)_\infty^3(q^9;q^9)_\infty^3}{(q^3;q^3)_\infty^7}+4q^2\frac{(q^6;q^6)_\infty^2(q^{18};q^{18})_\infty^3}{(q^3;q^3)_\infty^6}
    \end{align}
\end{lemma}
Equation (\ref{5}) is equivalent to the $3-$dissection of $\varphi(-q)$ \cite[Eqn. $14.3.2$]{hirschhorn2017powerofq}. Hirschhorn and Sellers \cite{hirschhornsellers} proved (\ref{6}).

\section{Congruences for $\overline{B}_{2^\alpha,3}(n)$, $\overline{B}_{2^{2\alpha},3^\beta}(n)$ and $\overline{B}_{2^{2\alpha +1},3^\beta}(n)$}

This section comprises the proof of theorems \ref{thm1},\ref{thm2} and \ref{thm3}.

\begin{proof}[Proof of Theorem \ref{thm1}]
    Substituting $\ell_1=2^\alpha$ and $\ell_2=3$ in (\ref{0}), we obtain 
    \begin{align}\label{th1e1}
        \sum_{n=0}^{\infty}\overline{B}_{2^\alpha,3}(n)q^n=\frac{(q^2;q^2)_\infty(q^{2^\alpha};q^{2^\alpha})_\infty^2(q^{3};q^{3})_\infty^2(q^{2^{\alpha+1}\cdot3};q^{2^{\alpha+1}\cdot 3})_\infty}{(q;q)_\infty^2(q^{2^{\alpha+1}};q^{2^{\alpha+1}})_\infty(q^{6};q^{6})_\infty(q^{2^\alpha\cdot 3};q^{2^\alpha\cdot 3})_\infty^2}. 
    \end{align}

    Using (\ref{1}), we will get 
    \begin{align*}
        \sum_{n=0}^{\infty}\overline{B}_{2^\alpha,3}(n)q^n=\frac{(q^4;q^4)_\infty^4(q^{12};q^{12})_\infty^2(q^{2^\alpha};q^{2^\alpha})_\infty^2(q^{2^{\alpha+1}\cdot 3};q^{2^{\alpha+1}\cdot 3})_\infty}{(q^2;q^2)_\infty^4(q^8;q^8)_\infty(q^{24};q^{24})_\infty(q^{2^{\alpha+1}};q^{2^{\alpha+1}})_\infty(q^{2^{\alpha}\cdot 3};q^{2^{\alpha}\cdot 3})_\infty^2}
         \end{align*}
        \begin{align}\label{th1e2}
             +2q \frac{(q^4;q^4)_\infty(q^6;q^6)_\infty(q^8;q^8)_\infty(q^{24};q^{24})_\infty(q^{2^\alpha};q^{2^\alpha})_\infty^2(q^{2^{\alpha+1}\cdot 3};q^{2^{\alpha+1}\cdot 3})_\infty}{(q^2;q^2)_\infty^3(q^{12};q^{12})_\infty(q^{2^{\alpha+1}};q^{2^{\alpha+1}})_\infty(q^{2^{\alpha}\cdot 3};q^{2^{\alpha}\cdot 3})_\infty^2}.
        \end{align}

        Extracting the terms involving $q^{2n}$ from both sides of (\ref{th1e2}) and then replacing $q^2$ by $q$ we get,
       \begin{align}\label{th1e3}
            \sum_{n=0}^{\infty}\overline{B}_{2^\alpha,3}(2n)q^n=\frac{(q^2;q^2)_\infty^4(q^6;q^6)_\infty^2(q^{2^{\alpha-1}};q^{2^{\alpha-1}})_\infty^2q^{2^{\alpha}\cdot 3};q^{2^{\alpha}\cdot 3})_\infty}{(q;q)_\infty^4 (q^4;q^4)_\infty (q^{12};q^{12})_\infty(q^{2^\alpha};q^{2^\alpha})_\infty(q^{2^{\alpha-1}\cdot 3};q^{2^{\alpha-1}\cdot 3})_\infty^2}.
       \end{align} 
       
        Employing (\ref{3}) in (\ref{th1e3}), we get

        \begin{align}\label{th1e4}
            \sum_{n=0}^{\infty}\overline{B}_{2^\alpha,3}(2n)q^n= \left( \frac{(q^4;q^4)_\infty^{14}}{(q^2;q^2)_\infty^{14}(q^8;q^8)_\infty^4}+4q\frac{(q^4;q^4)_\infty^2(q^8;q^8)_\infty^4}{(q^2;q^2)_\infty^{10}} \right) 
           \end{align}
           \begin{align*}
                \times \frac{(q^2;q^2)_\infty^4(q^6;q^6)_\infty^2(q^{2^{\alpha-1}};q^{2^{\alpha-1}})_\infty^2q^{2^{\alpha}\cdot 3};q^{2^{\alpha}\cdot 3})_\infty}{ (q^4;q^4)_\infty (q^{12};q^{12})_\infty(q^{2^\alpha};q^{2^\alpha})_\infty(q^{2^{\alpha-1}\cdot 3};q^{2^{\alpha-1}\cdot 3})_\infty^2}.
           \end{align*}
           
       After bringing out the terms containing $q^{2n+1}$ from both sides of (\ref{th1e4}), replace $q^2$ by $q$ to obtain,
       
       \begin{align}\label{th1e5}
          \sum_{n=0}^{\infty}\overline{B}_{2^\alpha,3}(4n+2)q^n =4\frac{(q^2;q^2)_\infty(q^4;q^4)_\infty^4(q^3;q^3)_\infty^2(q^{2^{\alpha-2}};q^{2^{\alpha-2}})_\infty^2 (q^{2^{\alpha-1}\cdot 3};q^{2^{\alpha-1}\cdot 3})_\infty}{(q;q)_\infty^6(q^6;q^6)_\infty (q^{2^{\alpha-1}};q^{2^{\alpha-1}})_\infty(q^{2^{\alpha-2}\cdot 3};q^{2^{\alpha-2}\cdot 3})_\infty^2}.
       \end{align}
       
   From (\ref{th1e5}) we will obtain (\ref{7}).
   Now, we substitute (\ref{1}) and (\ref{3}) in (\ref{th1e5})

 \begin{align}
 \label{th1e*1}      \sum_{n=0}^{\infty}\overline{B}_{2^\alpha,3}(4n+2)q^n =4\frac{(q^4;q^4)_\infty^{22}(q^{12};q^{12})_\infty^2(q^{2^{\alpha-2}};q^{2^{\alpha-2}})_\infty^2 (q^{2^{\alpha-1}\cdot 3};q^{2^{\alpha-1}\cdot 3})_\infty}
       {(q^2;q^2)_\infty^{18}(q^8;q^8)_\infty^5(q^{24};q^{24})_\infty(q^{2^{\alpha-1}};q^{2^{\alpha-1}})_\infty(q^{2^{\alpha-2}\cdot 3};q^{2^{\alpha-2}\cdot 3})_\infty^2}
       \end{align}
       \begin{align*}
            +8q \frac{(q^4;q^4)_\infty^{19}(q^6;q^6)_\infty(q^{24};q^{24})_\infty(q^{2^{\alpha-2}};q^{2^{\alpha-2}})_\infty^2 (q^{2^{\alpha-1}\cdot 3};q^{2^{\alpha-1}\cdot 3})}
       {(q^2;q^2)_\infty^{17}(q^{12};q^{12})_\infty(q^8;q^8)_\infty^3(q^{2^{\alpha-1}};q^{2^{\alpha-1}})_\infty(q^{2^{\alpha-2}\cdot 3};q^{2^{\alpha-2}\cdot 3})_\infty^2}
       \end{align*}
       \begin{align*}
            +16q \frac{(q^4;q^4)_\infty^{10}(q^{12};q^{12})_\infty^2(q^8;q^8)_\infty^3(q^{2^{\alpha-2}};q^{2^{\alpha-2}})_\infty^2 (q^{2^{\alpha-1}\cdot 3};q^{2^{\alpha-1}\cdot 3})}
       {(q^2;q^2)_\infty^{14}(q^{24};q^{24})_\infty(q^{2^{\alpha-1}};q^{2^{\alpha-1}})_\infty(q^{2^{\alpha-2}\cdot 3};q^{2^{\alpha-2}\cdot 3})_\infty^2}
       \end{align*}
       \begin{align*}
            +32q^2 \frac{(q^4;q^4)_\infty^{7}(q^6;q^6)_\infty(q^8;q^8)_\infty^5(q^{24};q^{24})_\infty(q^{2^{\alpha-2}};q^{2^{\alpha-2}})_\infty^2 (q^{2^{\alpha-1}\cdot 3};q^{2^{\alpha-1}\cdot 3})}
       {(q^2;q^2)_\infty^{13}(q^{12};q^{12})_\infty(q^{2^{\alpha-1}};q^{2^{\alpha-1}})_\infty(q^{2^{\alpha-2}\cdot 3};q^{2^{\alpha-2}\cdot 3})_\infty^2}.
       \end{align*}
   
Equate the odd powers on both sides of (\ref{th1e*1}) to obtain
              
       \begin{align}\label{th1e7}
       \sum_{n=0}^{\infty}\overline{B}_{2^\alpha,3}(8n+6)q^n 
            =8 \frac{(q^2;q^2)_\infty^{19}(q^3;q^3)_\infty(q^{12};q^{12})_\infty(q^{2^{\alpha-3}};q^{2^{\alpha-3}})_\infty^2 (q^{2^{\alpha-2}\cdot 3};q^{2^{\alpha-2}\cdot 3})}
       {(q;q)_\infty^{17}(q^{6};q^{6})_\infty(q^4;q^4)_\infty^3(q^{2^{\alpha-2}};q^{2^{\alpha-2}})_\infty(q^{2^{\alpha-3}\cdot 3};q^{2^{\alpha-3}\cdot 3})_\infty^2}
\end{align}
       \begin{align*}
            +16 \frac{(q^2;q^2)_\infty^{10}(q^{6};q^{6})_\infty^2(q^4;q^4)_\infty^3(q^{2^{\alpha-3}};q^{2^{\alpha-3}})_\infty^2 (q^{2^{\alpha-2}\cdot 3};q^{2^{\alpha-2}\cdot 3})}
       {(q;q)_\infty^{14}(q^{12};q^{12})_\infty(q^{2^{\alpha-2}};q^{2^{\alpha-2}})_\infty(q^{2^{\alpha-3}\cdot 3};q^{2^{\alpha-3}\cdot 3})_\infty^2}.
       \end{align*}
   Then (\ref{9}) follows from (\ref{th1e7}).        
       
      Isolating the terms having odd powers on both sides in (\ref{th1e2}) and then inserting (\ref{2}) we get,

        \begin{align}\label{th1e8}
            \sum_{n=0}^{\infty}\overline{B}_{2^\alpha,3}(2n+1)q^n=
            2 \left( \frac{(q^4;q^4)_\infty^6(q^6;q^6)_\infty^3}{(q^2;q^2)_\infty^9(q^{12};q^{12})_\infty^2}+3q\frac{(q^4;q^4)_\infty^2(q^6;q^6)_\infty(q^{12};q^{12})_\infty^2}{(q^2;q^2)_\infty^7}   \right)
           \end{align}

           \begin{align*}
           \times
             \frac{(q^2;q^2)_\infty(q^4;q^4)_\infty(q^{12};q^{12})_\infty(q^{2^{\alpha-1}};q^{2^{\alpha-1}})_\infty^2(q^{2^{\alpha}\cdot 3};q^{2^{\alpha}\cdot 3})_\infty}{(q^{6};q^{6})_\infty(q^{2^{\alpha}};q^{2^{\alpha}})_\infty(q^{2^{\alpha-1}\cdot 3};q^{2^{\alpha-1}\cdot 3})_\infty^2}.   
           \end{align*}
           
       Now we gather all the terms having even powers from both sides of (\ref{th1e8}) and combine it with (\ref{1}), (\ref{3}) and (\ref{4}). Then again by collecting the terms involving $q^{2n+1}$ and then replacing $q^2$ by $q$, one can obtain (\ref{8}).
  
\end{proof}

\begin{proof}[Proof of Theorem \ref{thm2}]

From both sides of (\ref{th1e8}) we identify and extract the terms containing $q^{2n+1}$ and then swap $q^2$ with $q$ to get

 \begin{align}\label{th2e3}
          \sum_{n=0}^{\infty}\overline{B}_{2^\alpha,3}(4n+3)q^n =6\frac{(q^2;q^2)_\infty^3(q^6;q^6)_\infty^3(q^{2^{\alpha-2}};q^{2^{\alpha-2}})_\infty^2 (q^{2^{\alpha-1}\cdot 3};q^{2^{\alpha-1}\cdot 3})_\infty}{(q;q)_\infty^6(q^{2^{\alpha-1}};q^{2^{\alpha-1}})_\infty(q^{2^{\alpha-2}\cdot 3};q^{2^{\alpha-2}\cdot 3})_\infty^2}.
       \end{align}

Consequently (\ref{th2e1}) follows for the case $k=0$.

Using (\ref{lp1}) in equation (\ref{th1e1}) gives
\begin{align}
        \sum_{n=0}^{\infty}\overline{B}_{2^\alpha,3}(n)q^n
        &\equiv \frac{(q^2;q^2)_\infty(q^{2^\alpha};q^{2^\alpha})_\infty^2(q;q)_\infty^6(q^{2^{\alpha+1}};q^{2^{\alpha+1}})_\infty^3}{(q;q)_\infty^2(q^{2^{\alpha+1}};q^{2^{\alpha+1}})_\infty(q^{2};q^{2})_\infty^3(q^{2^\alpha};q^{2^\alpha})_\infty^6}\\
 \label{th3e*3} &\equiv \frac{(q;q)_\infty^4 (q^{2^{\alpha+1}};q^{2^{\alpha+1}})_\infty^2}{(q^2;q^2)_\infty^2 (q^{2^\alpha};q^{2^\alpha})_\infty^4}\pmod 3.
    \end{align}
We can write (\ref{th3e*3}) as
\begin{align}
    \sum_{n=0}^{\infty}\overline{B}_{2^\alpha,3}(n)q^n
        &\equiv
        \left(\varphi(-q)\right)^2\frac{(q^{2^{\alpha+1}};q^{2^{\alpha+1}})_\infty^2}{(q^{2^\alpha};q^{2^\alpha})_\infty^4}\\
        &\equiv \left(1+2\sum_{i=1}^{\infty}(-q)^{i^2}\right)^2\frac{(q^{2^{\alpha+1}};q^{2^{\alpha+1}})_\infty^2}{(q^{2^\alpha};q^{2^\alpha})_\infty^4}\\
   \label{th3e*4}     &\equiv \left(1+\sum_{i=1}^{\infty}(-1)^{i^2}q^{i^2}+\sum_{i,j=1}^{\infty}(-1)^{i^2+j^2}q^{i^2+j^2}\right)\frac{(q^{2^{\alpha+1}};q^{2^{\alpha+1}})_\infty^2}{(q^{2^\alpha};q^{2^\alpha})_\infty^4}\pmod 3.
\end{align}
    Since $i^2\equiv 0\pmod 4$ if and only if $i$ is even, by extracting the terms containing $q^4$ from (\ref{th3e*4}) we can write,

    \begin{align}
    \sum_{n=0}^{\infty}\overline{B}_{2^\alpha,3}(4n)q^{4n}
      \label{th3e*+5}  &\equiv \left(\sum_{i=0}^{\infty}q^{4i^2}+\sum_{i,j=1}^{\infty}q^{4(i^2+j^2)}\right)\frac{(q^{2^{\alpha+1}};q^{2^{\alpha+1}})_\infty^2}{(q^{2^\alpha};q^{2^\alpha})_\infty^4}\pmod 3.
\end{align}
By replacing $q^4$ by $q$ in (\ref{th3e*+5}), we get
\begin{align}
     \sum_{n=0}^{\infty}\overline{B}_{2^\alpha,3}(4n)q^{n}
      \label{th3e*5}  &\equiv \left(\sum_{i=0}^{\infty}q^{i^2}+\sum_{i,j=1}^{\infty}q^{(i^2+j^2)}\right)\frac{(q^{2^{\alpha-1}};q^{2^{\alpha-1}})_\infty^2}{(q^{2^{\alpha-2}};q^{2^{\alpha-2}})_\infty^4}\pmod 3.
\end{align}
By using the facts that $i^2\equiv 0$ or $1\pmod 4$ and $i^2+j^2\equiv 0,1$ or $2\pmod 4$, we are able to conclude that there are no terms for the form $q^{4n+3}$ in (\ref{th3e*5}). Thus we can write if $\alpha>3$
\begin{align}
    \overline{B}_{2^\alpha,3}(4(4n+3))\equiv 0\pmod 3.
\end{align}
Since we have $\overline{B}_{2^\alpha,3}(4(4n+3))\equiv 0\pmod 2$, we will obtain (\ref{th2e1}) for the case $k=1$.

From both sides of (\ref{th3e*5}), if we extract the terms of the form $q^{4n}$, we will get

\begin{align}
    \sum_{n=0}^{\infty}\overline{B}_{2^\alpha,3}(4^2n)q^{4n}
    &\equiv \left(\sum_{i=0}^{\infty}q^{4i^2}+\sum_{i,j=1}^{\infty}q^{4(i^2+j^2)}\right)\frac{(q^{2^{\alpha-1}};q^{2^{\alpha-1}})_\infty^2}{(q^{2^{\alpha-2}};q^{2^{\alpha-2}})_\infty^4}\pmod 3
\end{align}

Replacing $q^4$ by $q$ we can write

\begin{align}
     \sum_{n=0}^{\infty}\overline{B}_{2^\alpha,3}(4^2n)q^{n}
    &\equiv \left(\sum_{i=0}^{\infty}q^{i^2}+\sum_{i,j=1}^{\infty}q^{i^2+j^2}\right)\frac{(q^{2^{\alpha-3}};q^{2^{\alpha-3}})_\infty^2}{(q^{2^{\alpha-4}};q^{2^{\alpha-4}})_\infty^4}\pmod 3.
\end{align}

Proceeding as above, if $\alpha>5$ we are able to write

\begin{align}
    \overline{B}_{2^\alpha,3}(4^2(4n+3))\equiv 0\pmod 6
\end{align}
which is the case for $k=2.$

If we continue to repeat this process, at the $k^{\text{th}}$ step we will obtain

\begin{align}
     \sum_{n=0}^{\infty}\overline{B}_{2^\alpha,3}(4^kn)q^{n}
    &\equiv \left(\sum_{i=0}^{\infty}q^{i^2}+\sum_{i,j=1}^{\infty}q^{i^2+j^2}\right)\frac{(q^{2^{\alpha+1-2k}};q^{2^{\alpha+1-2k}})_\infty^2}{(q^{2^{\alpha-2k}};q^{2^{\alpha-2k}})_\infty^4}\pmod 3.
\end{align}

Using the same arguments as above, if $\alpha>2k +1$, we get
\begin{align}
    \overline{B}_{2^\alpha,3}(4^k(4n+3))\equiv 0\pmod 6
\end{align}
which completes the proof of (\ref{th2e1})

Plugging (\ref{3}) and (\ref{4}) in (\ref{th2e3}) and interchanging $q^2$ with $q$ after collecting the terms comprising $q^{2n+1}$ on both sides, we get

\begin{align*}
    \sum_{n=0}^{\infty}\overline{B}_{2^\alpha,3}(8n+7)q^n =12\frac{(q^2;q^2)_\infty^{16}(q^8;q^8)_\infty^2(q^3;q^3)_\infty^3(q^{2^{\alpha-3}};q^{2^{\alpha-3}})_\infty^2(q^{2^{\alpha-2}\cdot 3};q^{2^{\alpha-2}\cdot 3})_\infty}{(q;q)_\infty^{16}(q^4;q^4)_\infty^5 (q^{2^{\alpha-2}};q^{2^{\alpha-2}})_\infty(q^{2^{\alpha-3}\cdot 3};q^{2^{\alpha-3}\cdot 3})_\infty^2}
\end{align*}

\begin{align*}
    +24 \frac{(q^4;q^4)_\infty^9(q^2;q^2)_\infty^2(q^3;q^3)_\infty^3 (q^{2^{\alpha-3}};q^{2^{\alpha-3}})_\infty^2(q^{2^{\alpha-2}\cdot 3};q^{2^{\alpha-2}\cdot 3})_\infty}{(q;q)_\infty^{12}(q^8;q^8)_\infty^2 (q^{2^{\alpha-2}};q^{2^{\alpha-2}})_\infty(q^{2^{\alpha-3}\cdot 3};q^{2^{\alpha-3}\cdot 3})_\infty^2}
\end{align*}

From this we obtain (\ref{th2e2}).

\end{proof}

\begin{proof}[Proof of Theorem \ref{thm3}]

Utilizing (\ref{0}) with $\ell_1=2^{2\alpha+1}$ and $\ell_2=3^{\beta}$, one can write

\begin{align}\label{th3e3}
    \sum_{n=0}^{\infty}\overline{B}_{2^{2\alpha +1},3^\beta}(n)q^n
    = \frac{(q^2;q^2)_\infty(q^{2^{2\alpha+1}};q^{2^{2\alpha+1}})_\infty^2(q^{3^\beta};q^{3^\beta})_\infty^2 (q^{2\cdot2^{2\alpha +1}\cdot 3^\beta};q^{2\cdot2^{2\alpha +1}\cdot 3^\beta})_\infty}
    {(q;q)_\infty^2 (q^{2\cdot2^{2\alpha +1}};q^{2\cdot2^{2\alpha +1}})_\infty (q^{2\cdot 3^\beta};q^{2\cdot 3^\beta})_\infty (q^{2^{2\alpha +1}\cdot 3^{\beta}};q^{2^{2\alpha +1}\cdot 3^{\beta}})_\infty^2}
\end{align}

Making use of (\ref{5}) with $q$ replaced with $q^{2^{2\alpha+1}}$ and (\ref{6}) we can write (\ref{th3e3}) as

\begin{align}\label{th3e4}
    \sum_{n=0}^{\infty}\overline{B}_{2^{2\alpha +1},3^\beta}(n)q^n
    = \left( \frac{(q^6;q^6)_\infty^4(q^9;q^9)_\infty^6}{(q^3;q^3)_\infty^8(q^{18};q^{18})_\infty^3}+2q\frac{(q^6;q^6)_\infty^3(q^9;q^9)_\infty^3}{(q^3;q^3)_\infty^7}+4q^2\frac{(q^6;q^6)_\infty^2(q^{18};q^{18})_\infty^3}{(q^3;q^3)_\infty^6} \right)
\end{align}
\begin{align*}
\times
     \left( \frac{(q^{9\cdot {2^{2\alpha+1}}};q^{9\cdot{2^{2\alpha+1}}})_\infty^2}{(q^{18\cdot {2^{2\alpha+1}}};q^{18\cdot {2^{2\alpha+1}}})_\infty}-2q^{2^{2\alpha+1}}\frac{(q^{3\cdot{2^{2\alpha+1}}};q^{3\cdot{2^{2\alpha+1}}})_\infty(q^{18\cdot{2^{2\alpha+1}}};q^{18\cdot{2^{2\alpha+1}}})_\infty^2}{(q^{6\cdot{2^{2\alpha+1}}};q^{6\cdot{2^{2\alpha+1}}})_\infty(q^{9\cdot{2^{2\alpha+1}}};q^{9\cdot{2^{2\alpha+1}}})_\infty} \right)
\end{align*}
   \begin{align*}
\times
\frac{(q^{3^\beta};q^{3^\beta})_\infty^2 (q^{2\cdot2^{2\alpha +1}\cdot 3^\beta};q^{2\cdot2^{2\alpha +1}\cdot 3^\beta})_\infty}
    {(q^{2\cdot 3^\beta};q^{2\cdot 3^\beta})_\infty (q^{2^{2\alpha +1}\cdot 3^{\beta}};q^{2^{2\alpha +1}\cdot 3^{\beta}})_\infty^2}
   \end{align*}

   Since $2^{2\alpha +1}+1\equiv 0\pmod 3$ and $2^{2\alpha +1}+2\not\equiv 0\pmod 3$, collecting the terms involving $q^{3n}$ from both sides of (\ref{th3e4}) and replacing $q^3$ by $q$,and then using (\ref{6}) we have

\begin{align}\label{th3e5}
     \sum_{n=0}^{\infty}\overline{B}_{2^{2\alpha +1},3^\beta}(3n)q^n
     \equiv   \frac{(q^6;q^6)_\infty^{13} (q^9;q^9)_\infty^{24} (q^{3\cdot 2^{2\alpha +1}};q^{3\cdot 2^{2\alpha +1}})_\infty^2 }
    {(q^3;q^3)_\infty^{26}(q^{18};q^{18})_\infty^{12} (q^{6\cdot 2^{2\alpha +1}};q^{6\cdot 2^{2\alpha +1}})_\infty}
\end{align}

\begin{align*}
    \times
    \frac{(q^{3^{\beta-1}};q^{3^{\beta-1}})_\infty^2 (q^{2\cdot2^{2\alpha +1}\cdot 3^{\beta-1}};q^{2\cdot2^{2\alpha +1}\cdot 3^{\beta-1}})_\infty}{(q^{2\cdot 3^{\beta-1}};q^{2\cdot 3^{\beta-1}})_\infty (q^{2^{2\alpha +1}\cdot 3^{\beta-1}};q^{2^{2\alpha +1}\cdot 3^{\beta-1}})_\infty^2} \pmod 4
\end{align*}
    From this we can see that (\ref{th3e1}) holds.

   Again by setting $\ell_1=2^{2\alpha}$, $\ell_2=3^\beta$ in (\ref{0}) we get 

   \begin{align}\label{th3e6}
    \sum_{n=0}^{\infty}\overline{B}_{2^{2\alpha },3^\beta}(n)q^n
    = \frac{(q^2;q^2)_\infty(q^{2^{2\alpha}};q^{2^{2\alpha}})_\infty^2(q^{3^\beta};q^{3^\beta})_\infty^2 (q^{2\cdot2^{2\alpha }\cdot 3^\beta};q^{2\cdot2^{2\alpha }\cdot 3^\beta})_\infty}
    {(q;q)_\infty^2 (q^{2\cdot2^{2\alpha }};q^{2\cdot2^{2\alpha }})_\infty (q^{2\cdot 3^\beta};q^{2\cdot 3^\beta})_\infty (q^{2^{2\alpha }\cdot 3^{\beta}};q^{2^{2\alpha }\cdot 3^{\beta}})_\infty^2}
\end{align}

Applying (\ref{5}) with $q$ replaced by $q^{2^{2\alpha}}$ and (\ref{6}) in (\ref{th3e6}), we will obtain

\begin{align}\label{th3e7}
    \sum_{n=0}^{\infty}\overline{B}_{2^{2\alpha },3^\beta}(n)q^n
    = \left( \frac{(q^6;q^6)_\infty^4(q^9;q^9)_\infty^6}{(q^3;q^3)_\infty^8(q^{18};q^{18})_\infty^3}+2q\frac{(q^6;q^6)_\infty^3(q^9;q^9)_\infty^3}{(q^3;q^3)_\infty^7}+4q^2\frac{(q^6;q^6)_\infty^2(q^{18};q^{18})_\infty^3}{(q^3;q^3)_\infty^6} \right)
\end{align}
\begin{align*}
\times
     \left( \frac{(q^{9\cdot {2^{2\alpha}}};q^{9\cdot{2^{2\alpha}}})_\infty^2}{(q^{18\cdot {2^{2\alpha}}};q^{18\cdot {2^{2\alpha}}})_\infty}-2q^{2^{2\alpha}}\frac{(q^{3\cdot{2^{2\alpha}}};q^{3\cdot{2^{2\alpha}}})_\infty(q^{18\cdot{2^{2\alpha}}};q^{18\cdot{2^{2\alpha}}})_\infty^2}{(q^{6\cdot{2^{2\alpha}}};q^{6\cdot{2^{2\alpha}}})_\infty(q^{9\cdot{2^{2\alpha}}};q^{9\cdot{2^{2\alpha}}})_\infty} \right)
\end{align*}
   \begin{align*}
\times
\frac{(q^{3^\beta};q^{3^\beta})_\infty^2 (q^{2\cdot2^{2\alpha }\cdot 3^\beta};q^{2\cdot2^{2\alpha }\cdot 3^\beta})_\infty}
    {(q^{2\cdot 3^\beta};q^{2\cdot 3^\beta})_\infty (q^{2^{2\alpha }\cdot 3^{\beta}};q^{2^{2\alpha }\cdot 3^{\beta}})_\infty^2}
   \end{align*}

   Because we have $2^{2\alpha }+2\equiv 0\pmod 3$ and $2^{2\alpha }+1\not\equiv 0\pmod 3$, taking all the terms containing $q^{3n}$ from sides of (\ref{th3e7}) and replacing $q^3$ by $q$, and then by using (\ref{6}) we can see that

\begin{align}\label{th3e8}
     \sum_{n=0}^{\infty}\overline{B}_{2^{2\alpha },3^\beta}(3n)q^n
     \equiv   \frac{(q^6;q^6)_\infty^{13} (q^9;q^9)_\infty^{24} (q^{3\cdot 2^{2\alpha }};q^{3\cdot 2^{2\alpha }})_\infty^2 }
    {(q^3;q^3)_\infty^{26}(q^{18};q^{18})_\infty^{12} (q^{6\cdot 2^{2\alpha }};q^{6\cdot 2^{2\alpha }})_\infty}
\end{align}

\begin{align*}
    \times
    \frac{(q^{3^{\beta-1}};q^{3^{\beta-1}})_\infty^2 (q^{2\cdot2^{2\alpha }\cdot 3^{\beta-1}};q^{2\cdot2^{2\alpha }\cdot 3^{\beta-1}})_\infty}{(q^{2\cdot 3^{\beta-1}};q^{2\cdot 3^{\beta-1}})_\infty (q^{2^{2\alpha }\cdot 3^{\beta-1}};q^{2^{2\alpha }\cdot 3^{\beta-1}})_\infty^2} \pmod 8
\end{align*}

Accordingly, (\ref{th3e8}) implies (\ref{th3e2}).\\

\end{proof}
 

\section{Congruences for $\overline{B}_{2^{\alpha},3^\beta}(n)$ }
This section is devoted for proving Theorem \ref{THM}.

\begin{proof}[Proof of Theorem \ref{THM}]
    Substituting $\ell_1=2^\alpha$ and $\ell_2=3^\beta$ in (\ref{0}), we can write
    \begin{align}\label{p1}
        \overline{B}_{2^\alpha, 3^\beta}(n)=\frac{(q^2;q^2)_\infty (q^{2^\alpha};q^{2^\alpha})_\infty^2(q^{3^\beta};q^{3^\beta})_\infty^2(q^{2^{\alpha+1}\cdot3^\beta};q^{2^{\alpha+1}\cdot3^\beta})_\infty}{(q;q)_\infty^2 (q^{2^{\alpha+1}};q^{2^{\alpha+1}})_\infty (q^{2\cdot3^\beta};q^{2\cdot3^\beta})_\infty (q^{2^{\alpha}\cdot3^\beta};q^{2^{\alpha}\cdot3^\beta})_\infty^2}
    \end{align}

     Since $\varphi(-q)=\frac{(q;q)_\infty^2}{(q^2;q^2)_\infty}$, we can rewrtite \ref{p1} as 
    \begin{align}\label{p2}
        \sum_{n=0}^{\infty}\overline{B}_{2^{\alpha },3^\beta}(n)q^n=\frac{\varphi(-q^{2^\alpha})\varphi(-q^{3^\beta})}{\varphi(-q)\varphi(-q^{2^{\alpha}\cdot3^{\beta}})}.
    \end{align}
    We use the facts that $\frac{1}{\varphi(-q)}=\varphi(q)(\varphi(q^2))^2(\varphi(q^4))^4(\varphi(q^8))^8\cdots$ and $(\varphi(q^i))^j\equiv 1\pmod 8$ if $j\geq 4$ and $j$ is a multiple of $4$, to rewrite (\ref{p2}) as

    \begin{align}\label{p3}
        \sum_{n=0}^{\infty}\overline{B}_{2^{\alpha },3^\beta}(n)q^n\equiv \frac{\varphi(-q^{2^\alpha})\varphi(-q^{3^\beta})\varphi(q)(\varphi(q^2))^2}{\varphi(-q^{2^{\alpha}\cdot3^{\beta}})}\pmod 8.
    \end{align}
For $\alpha\geq 2,$ $\varphi(-q^{2^{\alpha}\cdot3^{\alpha}})$ contains only the terms of the form $q^{12k}$. Now we consider the numerator part of (\ref{p3})

\begin{align}
   \notag &\varphi(-q^{2^\alpha})\varphi(-q^{3^\beta})\varphi(q)(\varphi(q^2))^2\\
  \notag & =\left(1+2\sum_{n=1}^{\infty}(-q^{2^\alpha})^{n^2} \right) \left(1+2\sum_{n=1}^{\infty}(-q^{3^\beta})^{n^2} \right) \left(1+2\sum_{n=1}^{\infty}q^{n^2} \right)
    \left(1+2\sum_{n=1}^{\infty}(q^{2})^{n^2} \right)^2\\
 \notag  & \equiv 1+ 2\sum_{n=1}^{\infty}(-q^{2^\alpha})^{n^2} +2 \sum_{n=1}^{\infty}(-q^{3^\beta})^{n^2} + 2 \sum_{n=1}^{\infty}q^{n^2}
   + 4 \sum_{n=1}^{\infty}(q^{2})^{n^2} + 4 \left(\sum_{n=1}^{\infty}(q^{2})^{n^2}\right)^2\\
\label{p4} & + 4 \sum_{i,j=1}^{\infty}(-q^{2^\alpha})^{i^2}(-q^{3^\beta})^{j^2} 
    + 4 \sum_{i,j=1}^{\infty}(-q^{2^\alpha})^{i^2}q^{j^2}
    + 4 \sum_{i,j=1}^{\infty}(-q^{3^\beta})^{i^2}q^{j^2} \pmod 8.
\end{align}

Since we have the facts that, 
  $  n^2 \equiv 0,1,4 \text{ or } 9 \pmod{12}, 
    2^\alpha \equiv 4 \text{ or } 8 \pmod{12} \text{ and }
    3^\beta \equiv 3 \text{ or } 9 \pmod{12}$, we can conclude the following.

If $\alpha$ is odd and $\beta$ is even,
\begin{align*}
  2^\alpha  n^2 &\equiv 0 \text{ or } 8 \pmod{12}\\ 
  3^\beta n^2 &\equiv 0 \text{ or } 9 \pmod{12}\\
    2^\alpha i^2 +j^2 &\equiv0,1 ,4,5,8 \text{ or } 9 \pmod{12} \\
    3^\beta i^2 +j^2 &\equiv 0,1 ,4,6,9 \text{ or } 10 \pmod{12}\\
  2^\alpha i^2  +3^\beta j^2 &\equiv 0,5,8 \text{ or } 9 \pmod{12}.
\end{align*}

When $\alpha$ is even and $\beta$ is odd,
\begin{align*}
  2^\alpha  n^2 &\equiv 0 \text{ or } 4 \pmod{12}\\ 
  3^\beta n^2 &\equiv 0 \text{ or } 3 \pmod{12}\\
    2^\alpha i^2 +j^2 &\equiv0,1 ,4,5,8 \text{ or } 9 \pmod{12} \\
    3^\beta i^2 +j^2 &\equiv 0,1 ,3,4,7 \text{ or } 9 \pmod{12}\\
  2^\alpha i^2  +3^\beta j^2 &\equiv 0,3,4 \text{ or } 7 \pmod{12}.
\end{align*}

If $\alpha$ and $\beta$ are odd,
\begin{align*}
  2^\alpha  n^2 &\equiv 0 \text{ or } 8 \pmod{12}\\ 
  3^\beta n^2 &\equiv 0 \text{ or } 3 \pmod{12}\\
    2^\alpha i^2 +j^2 &\equiv0,1 ,4,5,8 \text{ or } 9 \pmod{12} \\
    3^\beta i^2 +j^2 &\equiv  0,1 ,3,4,7\text{ or } 9 \pmod{12}\\
  2^\alpha i^2  +3^\beta j^2 &\equiv 0,3,8 \text{ or } 11 \pmod{12}.
\end{align*}

If $\alpha$ and $\beta$ are even,
\begin{align*}
  2^\alpha  n^2 &\equiv 0 \text{ or } 4 \pmod{12}\\ 
  3^\beta n^2 &\equiv 0 \text{ or } 9 \pmod{12}\\
    2^\alpha i^2 +j^2 &\equiv0,1 ,4,5,8 \text{ or } 9 \pmod{12} \\
    3^\beta i^2 +j^2 &\equiv 0,1,4,6,9 \text{ or } 10 \pmod{12}\\
  2^\alpha i^2  +3^\beta j^2 &\equiv 0,1,4 \text{ or } 9 \pmod{12}.
\end{align*}
  
If $\alpha$ is odd and $\beta$ is even, in (\ref{p4}) there exist no terms of the form $q^{12n+3}$, $q^{12n+7}$ and $q^{12n+11}$. When $\alpha$ is even and $\beta$ is odd, no terms of the form $q^{12n+11}$ can be found. For the case where both $\alpha$ and $\beta$ are even, there does not present any terms of the form $q^{12n+3}$, $q^{12n+7}$ and $q^{12n+11}$. From these discussions, we can deduce Theorem \ref{THM}.
    
\end{proof}

\end{document}